\RequirePackage{fix-cm}
\documentclass[smallextended]{svjour3}       
\smartqed  
\usepackage{graphicx}
 \journalname{Periodica Mathematica Hungarica}
\begin{document}

\title{Approximate solution of a system of singular integral
equations of the first kind by using Chebyshev polynomials
}

\titlerunning{Approximate solution of a system of singular integral
equations}        

\author{S. Shahmorad \and S. Ahdiaghdam 
}


\institute{S. Shahmorad \at
              Faculty of Mathematical
Sciences, University of Tabriz, Tabriz, Iran \\
              Tel.: +98-914-302-6376\\
              Fax:  +98-413-334-2102\\
              \email{shahmorad@tabrizu.ac.ir}           
           \and
           S. Ahdiaghdam \at
              Islamic Azad University - Marand Branch \\
              \email{ahdi@marandiau.ac.ir}
}

\date{Received: date / Accepted: date}

\maketitle

\begin{abstract}
The aim of the present work is to introduce a method based on
Chebyshev polynomials for the numerical solution of a system of
Cauchy type singular integral equations of the first kind on a
finite segment. Moreover, an estimation error is computed for the
approximate solution. Numerical results demonstrate effectiveness of
the proposed method. 
\end{abstract}

\section{Introduction}
\label{intro} Let us consider a system of singular integral
equations of the form
\begin{eqnarray}\label{1}
A(t) \Phi(t) + \int_{-1}^{1} \frac{B(\tau)
\Phi(\tau)}{\tau-t}\,d\tau +  \int_{-1}^{1} K(t,\tau) \Phi(\tau)
d\tau =F(t)\,, & -1<t<1
\end{eqnarray}
where
\begin{eqnarray*} K(t,\tau)=\left[ K_{ij}(t,\tau)\right]
,\,\,\,\,i,j=1,2,\ldots,N\,,\\
F(t)=[f_1(t)\,,\,f_2(t)\,,\,\ldots\,,\,f_N(t)]^T\,,\\
\Phi(t)=[\phi_1(t)\,,\,\phi_2(t)\,,\,\ldots\,,\,\phi_N(t)]^T\,,\\
A(t)=\left[ a_{ij}(t)\right] \,,\,\,\,\,i,j=1,2,\ldots,N\,, \\
B(t)=\left[ b_{ij}(t)\right] \,,\,\,\,\,i,j=1,2,\ldots,N\,.
\end{eqnarray*}
Here, $\{K_{ij}(t,\tau)\}_{i,j=1}^N$ and $\{f_i(t)\}_{i=1}^N$ are
given real-valued H$\ddot{o}$lder functions and
$\{\phi_j(t)\}_{j=1}^N$ are the unknown functions. The matrices $A$
and $B$ are known such that $S=A+B$ and $D=A-B$ are nonsingular for
all $t\in[-1, 1]$. In some familiar physical problems the entries of
the matrices $A$ and $B$ are constants.

The singular integral equations play important roles in physics and
theoretical mechanics, particularly in the areas of elasticity,
aerodynamics, and unsteady aerofoil theory. They are highly
effective in solving boundary value problems occurring in the theory
of functions of a complex variable, potential theory, the theory of
elasticity, and the theory of fluid mechanics. A general theory of
the system of equations (\ref{1}) has given in \cite{mu}.

We study the system (\ref{1}) in the case that $A(t)=0$ and $B(t)$
is a constant matrix. Therefore, the $i^{th}$ equation of system
(\ref{1}) takes the form
\begin{eqnarray}\label{2}
\int_{-1}^{1} \sum_{j=1}^N b_{ij} \phi_j(\tau)
\frac{d\tau}{\tau-t}\, +  \int_{-1}^{1} \sum_{j=1}^N K_{ij}(t,\tau)
\phi_j(\tau) d\tau =f_i(t)\,, \ \ \ -1<t<1\,.
\end{eqnarray}

Studies on this singular integral equation can be found in some
literatures (See \cite{ab,cb,kp,ks}). Chakrabarti and Berghe
\cite{cb} proposed a method for solving Eq. (\ref{2}) using
polynomial approximation and collocation points have chosen to be
the zeros of the Chebyshev polynomials of the first kind for all
cases. Kashfi and Shahmorad \cite{ks} constructed another
approximate solution of this equation using Chebyshev polynomials of
the first and second kinds. Some other methods for solving this
equation can be found in \cite{ab,kp}. A convergence analysis of
Galerkin and collocation methods for Eq. (\ref{2}) has been given by
Miel \cite{mi}.

A special type of Eq. (\ref{2}) is the famous Cauchy singular
integral equation
\begin{eqnarray}\label{3}
\int_{-1}^{1} \frac{\phi(\tau)}{\tau-t}\,d\tau =f(t)\,, \ \ \
 -1<t<1\,,
\end{eqnarray}
which has the following analytical solutions in four special cases
based on boundedness of the unknown function $ \phi(\tau) $ at the
endpoints of the interval $[-1,\,1]$ \cite{cb,mc,se}. \\ \\
{\bf Case 1.} If the function $ \phi(\tau) $ is unbounded at the
 endpoints $ \tau=\pm1 $, then
\begin{eqnarray}\label{4}
\phi(\tau)=\frac{a_0}{\sqrt{1-\tau^2}} - \frac{1}{{\pi}^2
\sqrt{1-\tau^2}} \int_{-1}^{1} \frac{\sqrt{1-t^2}
f(t)}{t-\tau}\,dt\,, \ \ \ -1<\tau<1
\end{eqnarray}
where $a_0$ is an arbitrary constant. \\ \\
{\bf Case 2.} If the function $ \phi(\tau) $ is bounded at the
endpoints $ \tau=\pm1 $, then
\begin{eqnarray}\label{5}
\phi(\tau)= - \frac{\sqrt{1-\tau^2}}{{\pi}^2 } \int_{-1}^{1} \frac{
f(t)}{\sqrt{1-t^2}\,(t-\tau)}\,dt\,, \ \ \ -1<\tau<1
\end{eqnarray}
a necessary and sufficient condition of existing this solution is:
\begin{equation}\label{6}
\int_{-1}^{1} \frac{f(t)}{\sqrt{1-t^2}}\,dt =0\,.
\end{equation} \\ \\
{\bf Case 3.} If the function $ \phi(\tau) $ is bounded at the
endpoint $ \tau=-1 $ and unbounded at the endpoint $ \tau=1 $, then
\begin{eqnarray}\label{7}
\phi(\tau)=-\frac{1}{{\pi}^2}\,\sqrt{\frac{1+\tau}{1-\tau}}\,
\int_{-1}^{1} \sqrt{\frac{1-t}{1+t}}\frac{ f(t)}{t-\tau}\,dt\,, &
-1<\tau<1 .
\end{eqnarray}\\ \\
{\bf Case 4.} If the function $ \phi(\tau) $ is bounded at the
endpoint $ \tau=1 $ and unbounded at the endpoint $ \tau=-1 $, then
\begin{eqnarray}\label{8}
\phi(\tau)=-\frac{1}{{\pi}^2}\,\sqrt{\frac{1-\tau}{1+\tau}}\,
\int_{-1}^{1} \sqrt{\frac{1+t}{1-t}}\frac{ f(t)}{t-\tau}\,dt\,, &
-1<\tau<1 .
\end{eqnarray}
More methods for solving Eq. (\ref{3}) have given in
\cite{en,lz,se,ts}.

In the next section, we investigate approximate solutions for system
(\ref{1}) in the above four cases.
\section{Approximate solution}
To find approximate solutions for system (\ref{1}) in cases
$\textbf{1,2,3,4}$, for $\nu\in\{1,2,3,4\}$ we set
\begin{equation}\label{9}
\phi_j(\tau)\simeq
\varphi_{\nu,j}(\tau):=\frac{\lambda_\nu(\tau)}{\sqrt{1-{\tau}^2}}\sum_{l=0}^{M}
\beta_{jl}\,P_{\nu,l}(\tau)\,, \ \ \ j=1,2,\ldots ,N
\end{equation}
and
\begin{equation}\label{10}
K_{ij}(t,\tau):=\sum_{k=0}^{M} \gamma_{ijk}(t)\,P_{\nu,k}(\tau)\,, \
\ \ \ i,j=1,2,\ldots ,N\,,
\end{equation}
where
\begin{eqnarray*}
P_{\nu,j}(x)=\left\lbrace \begin{array}{cc}
T_j(x) =\cos(j \theta),\ \ &\nu=1,\\ \\
U_{j}(x) =\frac{\sin\left((j+1)\theta\right)
}{\sin(\theta)},\ \ &\nu=2,\\ \\
V_j(x) =\frac{\cos\left((j+\frac{1}{2})\theta\right)
}{\cos(\frac{\theta}{2})},\ \ &\nu=3,\\ \\
W_j(x) =\frac{\sin\left((j+\frac{1}{2})\theta\right)
}{\sin(\frac{\theta}{2})},\ \ &\nu=4,
\end{array}\right.
\end{eqnarray*}
are the Chebyshev polynomials of the first to fourth kinds and
\begin{eqnarray*}
\lambda_\nu(t)=\left\lbrace \begin{array}{cc}
 1,\ \ \ \ \ \ \ &\nu=1,\\
1-t^2,\ \ &\nu=2,\\ 1+t,\ \ &\nu=3,\\ 1-t,\ \ &\nu=4,
\end{array}
\right.
\end{eqnarray*}
in which $x=\cos(\theta)$. The roots of Chebyshev polynomials
$P_{\nu,M+1}(x)$ are given by
\begin{eqnarray}\label{11}
x_{\nu,n}=\left\lbrace \begin{array}{cc}\cos\left(
\frac{(2n-1)\pi}{2(M+1)}\right), \ \ &\nu=1,\\ \\
\cos\left(\frac{n\pi}{M+2}\right), \ \ &\nu=2,\\ \\
\cos\left(\frac{(2n-1)\pi}{2M+3}\right), \ \ &\nu=3,\\ \\
\cos\left(\frac{2n\pi}{2M+3}\right), \ \ &\nu=4,
\end{array}\right.
\end{eqnarray}
where $n=1,2,\ldots\,,M+1$. These roots are used as the nods of
Gauss-Chebyshev quadrature rules.
\\
\textbf{Lemma.}\cite{mh} The Chebyshev polynomials satisfy the
orthogonality property
\begin{equation}\label{12}
\int_{-1}^{1} \frac{\lambda_\nu(t)}{\sqrt{1-t^2}} P_{\nu,i}(t)
P_{\nu,j}(t)\,dt\, =\left\{\begin{array}{lll}
        0,\ \ & i \neq j, \\ \\
         \pi,\ \  & i=j=0\,, \ \ &\nu=1,\\ \\
          \frac{\pi}{2},\ \  & i=j\neq 0,\ \ &\nu=1,\\ \\
       \frac{\pi}{2},\ \  & i=j,\ \ &\nu=2,\\ \\
          {\pi},\ \  & i=j,\ \ &\nu=3\,,\,4 .
       \end{array}\right.
\end{equation}
\\
\textbf{Theorem 1.}\cite{mh} As a Cauchy principle value integral,
we have
\begin{eqnarray}\label{13}
\int_{-1}^{1}\frac{\lambda_\nu(\tau)
}{\sqrt{1-\tau^2}}\,\frac{P_{\nu,j}(\tau)}{\tau-t}\,d\tau\, =
{\pi}\, \left\{\begin{array}{cc}
 U_{j-1}(t), \ \ & \nu=1,\\ \\
-T_{j+1}(t), \ \ & \nu=2,\\  \\
 W_{j}(t), \ \ & \nu=3, \\ \\
 -V_{j}(t), \ \ & \nu=4.
\end{array}\right.
\end{eqnarray}
Now we describe details of finding approximate solution in cases
$\textbf{1-4}$.

{\bf Case 1.} For $\nu=1$ the relations (\ref{9})-(\ref{10}) take
the forms
\begin{equation}\label{a1}
\phi_j(\tau)\simeq
\varphi_{1,j}(\tau):=\frac{1}{\sqrt{1-{\tau}^2}}\sum_{l=0}^{M}{}^{'}
\beta_{jl}\,T_l(\tau)\,, \ \ \ \ j=1,2,\ldots ,N
\end{equation}
and
\begin{equation}\label{a2}
K_{ij}(t,\tau):=\sum_{k=0}^{M}{}^{'} \gamma_{ijk}(t)\,T_k(\tau)\,,\
\ \ \ \ i,j=1,2,\ldots ,N
\end{equation}
where $ \beta_{jl} \ \ (j=1,2,\ldots ,N,\ \ \ l=0,1,\ldots ,M)$ are
unknown coefficients and the symbol ($\sum {}^{'}$) denotes that the
first term in the summation is halved. The functions
$$\gamma_{ijk}(t)=\frac{2}{\pi}\int_{-1}^{1}\frac{K_{ij}(t,\tau) T_{k}(\tau)}{\sqrt{1-{\tau}2}}\,d\tau\,, \ \ \ \
i,j=1,2,\ldots ,N,\ \ \ \ k=0,1,\ldots ,M$$ can be determined
exactly, or may be approximated by using Gauss-Chebyshev quadrature
rule, i.e
$$\gamma_{ijk}(t)\simeq\frac{2}{M+1}\sum_{s=1}^{M+1}  K_{ij}(t,x_{1,s}) T_k(x_{1,s})\,,
$$
where $ x_{1,s}$ obtain from (\ref{11}).

Substituting from (\ref{a1})-(\ref{a2}) into Eq. (\ref{2}) and using
(\ref{12})-(\ref{13}) for $\nu=1$, gives the system
\begin{equation}\label{a3}
\sum_{j=1}^{N}\sum_{l=1}^{M}  b_{ij} \beta_{jl} U_{l-1}(t) +
\frac{1}{2} \sum_{j=1}^{N}\sum_{k=0}^{M}{}^{'} \gamma_{ijk}(t)
\beta_{jk} = \frac{1}{\pi} f_i(t)\,,\ \ \ \ \ i=1,2,\ldots ,N .
\end{equation}

If the given functions $ f_i(t) $ and  $ \gamma_{ijk}(t) $ are
square integrable on $ [-1\,,\,\,1] $ with respect to the weight
function $ \frac{\lambda_1(t)}{\sqrt{1-t^2}} $, then they can be
expanded as
\begin{equation}\label{a4}
\left\{\begin{array}{lll} \gamma_{ijk}(t)\simeq\sum_{l=0}^{M-1}
G_{ijkl}\,U_{l}(t) \,,\ \ \ \ \ \ &i,j=1,2,\ldots ,N,&\ \ \ \ \
k=0\,,1\,,\ldots\,,M \\ \\ \frac{1}{\pi}
f_i(t)\simeq\sum_{l=0}^{M-1} c_{il}\,U_{l}(t) \,,\ \ \ \ \
&i=1,2,\ldots ,N,&
\end{array}
\right.
\end{equation}
where the coefficients
\begin{eqnarray*}
\left\{\begin{array}{lll}
G_{ijkl}&=\frac{2}{\pi}\int_{-1}^{1}\sqrt{1-t^2}\gamma_{ijk}(t)U_{l}(t)\,dt\\
 &=\frac{4}{{\pi}^2}\int_{-1}^{1} \int_{-1}^{1}
\sqrt{\frac{1-t^2}{1-{\tau}^2}}K_{ij}(t,\tau)
U_{l}(t)\,T_{k}(\tau)\, d\tau dt& \\
&  \ \ \ \ \ \ \ \ \ \ \ \
 i,j=1,2,\ldots ,N, \ \ \ k=0,1,\ldots,M\,, \ \ &l=0,1,\ldots,M-1\, \\ \\
\ c_{il}&=\frac{1}{{\pi}^2}\int_{-1}^{1}\sqrt{1-t^2}f_i(t)
U_{l}(t)\,dt\,,\ \  i=1,2,\ldots ,N, \ \ &l=0,1,\ldots ,M-1,
 \end{array}\right.
\end{eqnarray*}
 can be approximately determined from
 \begin{equation}\label{a5}
\left\{\begin{array}{ll} G_{ijkl}\simeq\frac{4}{
(M+1)^2}\sum_{r=1}^{M} \sum_{s=1}^{M+1}(1-x_{2,r}^2)
K_{ij}(x_{2,r},x_{1,s}) U_{l}(x_{2,r}) T_k(x_{1,s})\,, \\ \\
c_{il}\simeq\frac{2}{\pi (M+1)}\sum_{r=1}^{M}(1-x_{2,r}^2)
f_i(x_{2,r}) U_{l}(x_{2,r}) .\end{array}\right.
\end{equation}
Using (\ref{a4}) in (\ref{a3}) and linearly independence of
$\{U_l(t)\}$, yield
$$\sum_{l=0}^{M-1}\sum_{j=1}^{N} b_{ij} \beta_{j\{l+1\}} U_{l}(t)  + \frac{1}{2}
\sum_{l=0}^{M-1}\sum_{j=1}^{N}\sum_{k=0}^{M}{}^{'} G_{ijkl}
\beta_{jk} U_{l}(t)= \sum_{l=0}^{M-1}  c_{il}\,U_{l}(t) ,$$ which
leads to the linear system
\begin{equation}\label{a6}
\sum_{j=1}^{N} \left[ b_{ij} \beta_{j\{l+1\}}
 + \frac{1}{2}
\sum_{k=0}^{M}{}^{'} G_{ijkl} \beta_{jk}\right] = c_{il}, \ \
i=1,2,\ldots ,N, \ \ l=0,1,\ldots ,M-1
\end{equation}
for the unknown values $ \beta_{jk} \ \ (j=1,2,\ldots ,N,\ \ \
k=0,1,\ldots ,M)$.  By taking arbitrary values for example for
$\beta_{11},\ldots,\beta_{N1}$, the remaining coefficients
$\beta_{jk}$ are uniquely found via the linear system (\ref{a6})
which determine the elements of the vector
function $\Phi(t)$ via Eq. (\ref{a1}).\\

{\bf Case 2.} We set $\nu=2$ in (\ref{9})-(\ref{10}) and substitute
them in Eq. (\ref{2}) to get
\begin{equation}\label{b3}
-\sum_{j=1}^{N}\sum_{l=0}^{M} b_{ij} \beta_{jl} T_{l+1}(t) +
\frac{1}{2} \sum_{j=1}^{N}\sum_{k=0}^{M}  \gamma_{ijk}(t) \beta_{jk}
= \frac{1}{\pi} f_i(t)\,,\ \ \ \ \ i=1,2,\ldots ,N .
\end{equation}
where we used the formulae (\ref{12})-(\ref{13}). Then, we expand
the functions $ f_i(t) $ and $ \gamma_{ijk}(t) $ as
\begin{eqnarray*}
\left\{\begin{array}{lll} \gamma_{ijk}(t)\simeq\sum_{l=0}^{M}{}^{'}
G_{ijkl}\,T_l(t)\,, \ \ \ &i,j=1,2,\ldots ,N, \ \ &
k=0\,,1\,,\ldots\,,M \\ \\
\frac{1}{\pi} f_i(t)\simeq\sum_{l=0}^{M}{}^{'} c_{il}\,T_l(t)\,,\ \
\ & \ i=1,2,\ldots ,N \ \ & \end{array} \right.
\end{eqnarray*}
where the coefficients are determined by
\begin{eqnarray*}
\left\{\begin{array}{ll}
G_{ijkl}=\frac{2}{\pi}\int_{-1}^{1}\frac{1}{\sqrt{1-t^2}}\gamma_{ijk}(t)T_{l}(t)\,dt\,,
\ \ \ \ \  i,j=1,2,\ldots ,N,\ \ \ \ \ &k,l=0\,,1\,,\ldots\,,M\,, \\
\ \ \ \ \ \ \  =\frac{4}{{\pi}^2}\int_{-1}^{1} \int_{-1}^{1}
\sqrt{\frac{1-\tau^2}{1-t^2}}
 K_{ij}(t,\tau) T_l(t)\,U_{k}(\tau)\, d\tau dt & \\ \\
c_{il}=\frac{2}{{\pi}^2}\int_{-1}^{1}\frac{1}{\sqrt{1-t^2}}f_i(t)
T_l(t)\,dt\,,\ \ \ \ \ \ \ \ \ \ i=1,2,\ldots ,N,\ \ \
&l=0\,,1\,,\ldots\,,M
 \end{array}\right.
 \end{eqnarray*}
 or approximated by
\begin{equation}\label{b4}
\left\{\begin{array}{ll}
G_{ijkl}\simeq\frac{4}{(M+1)(M+2)}\sum_{r=1}^{M+1}
\sum_{s=1}^{M+1}(1-x_{2,s}^2) K_{ij}(x_{1,r},x_{2,s}) T_l(x_{1,r})
U_k(x_{2,s}) \,, \\ \\
c_{il}\simeq\frac{2}{\pi (M+1)}\sum_{r=1}^{M+1} f_i(x_{1,r})
T_l(x_{1,r}).\end{array}\right.
\end{equation}
Using the last expansions in Eq. (\ref{b3}), returns the following
linear system of equations
\begin{eqnarray}\label{b6}
\left\lbrace \begin{array}{lll} \frac{1}{2} \sum_{j=1}^{N}
\sum_{k=0}^{M} G_{ijkl} \beta_{jk} =  c_{il}, \ \ &
i=1,2,\ldots\,,N, \ \ & l=0
\\ \\ \sum_{j=1}^{N} \left\lbrace - b_{ij} \beta_{j\{l-1\}} + \frac{1}{2}
\sum_{k=0}^{M} G_{ijkl} \beta_{jk}\right\rbrace = c_{il}, \ \  &
i=1,2,\ldots\,,N, \ \ & l=1,2,\ldots ,M
\end{array} \right.
\end{eqnarray}
for the unknown values $ \beta_{jl} \ \ (j=1,2,\ldots ,N,\ \ \
l=0,1,\ldots ,M)$. Then the elements of the vector function
$\Phi(t)$ obtain from Eq. (\ref{9}).
\\
{\bf Case 3,4.} Proceeding by the same way as we did in cases $1,2$,
we get the linear systems
\begin{eqnarray}\label{c6}
 \sum_{j=1}^{N} \left\lbrace  b_{ij} \beta_{jl} + \sum_{k=0}^{M}
 G_{ijkl} \beta_{jk}\right\rbrace =  c_{il}, \ \ \   i=1,2,\ldots\,,N,\ \ \ l=0,1,\ldots
 ,M
\end{eqnarray}
and
\begin{eqnarray}\label{d6}
\sum_{j=1}^{N} \left\lbrace - b_{ij} \beta_{jl} + \sum_{k=0}^{M}
G_{ijkl} \beta_{jk}\right\rbrace =  c_{il}, \ \ \ i=1,2,\ldots\,,N,\
\ \ l=0,1,\ldots ,M
\end{eqnarray}
respectively for $\nu=3$ and $\nu=4$, and then we determine the
elements of corresponding vector $\Phi(t)$ via (\ref{9}).

\section{An estimation error and numerical results}
In this section, we describe an estimation error for the approximate
solution. Let
$$\overline{\Phi}(t)=[\varphi_1(t)\,,\,\varphi_2(t)\,,\,\ldots\,,\,\varphi_N(t)]^T
$$ be the vector of approximate solution of the system (\ref{1}) and
$E(t)=\overline{\Phi}(t)-\Phi(t) $  be the associated vector valued
error function. Due to the approximation $\overline{\Phi}(t)$, for
$A(t)=0$ the system (\ref{1}) may be written as
\begin{equation}\label{E1} \int_{-1}^{1}
\frac{B(\tau) \overline{\Phi}(\tau)}{\tau-t}\,d\tau + \int_{-1}^{1}
K(t,\tau) \overline{\Phi}(\tau) d\tau =F(t)+H(t)\,, \ \ \ \ -1<t<1 ,
\end{equation}
where the perturbation term $H(t)$ obtains from
\begin{equation}\label{E2}
H(t)=\int_{-1}^{1} \frac{B(\tau)
\overline{\Phi}(\tau)}{\tau-t}\,d\tau
 +  \int_{-1}^{1} K(t,\tau) \overline{\Phi}(\tau) d\tau - F(t)\,, \
\ \ \ \ -1<t<1 .
\end{equation}
Subtracting Eq. (\ref{E1}) from Eq. (\ref{1}), yields a system of
error equations as
\begin{equation}\label{E3}
\int_{-1}^{1} \frac{B(\tau) {E(\tau)}}{\tau-t}\,d\tau +
\int_{-1}^{1} K(t,\tau) {E(\tau)} d\tau =H(t)\,, \ \ \ \  -1<t<1
\end{equation}
which is solvable approximately like the system (\ref{1}).\\

The following examples illustrate application of the method.\\

\textbf{Example $1.$} Let
\begin{eqnarray*}
A(t)=0,\ \ \ B(t)=\left(
        \begin{array}{cc}
        1 &\ \ 0 \\
        0 &\ \ 1 \\
        \end{array}
        \right), \ \ \
K(t,\tau)=\left(
        \begin{array}{cc}
        \tau-t &\ \ t \\
        \tau &\ \ \tau+t \\
        \end{array}
        \right),\ \ \ \ f_1(t)=\pi, \ \ \
        f_2(t)=2\pi t
\end{eqnarray*}
and find the solution of system (\ref{1}) in case $\textbf{1}$.

\begin{figure}
\begin{center}
\includegraphics[width=0.4\textwidth]{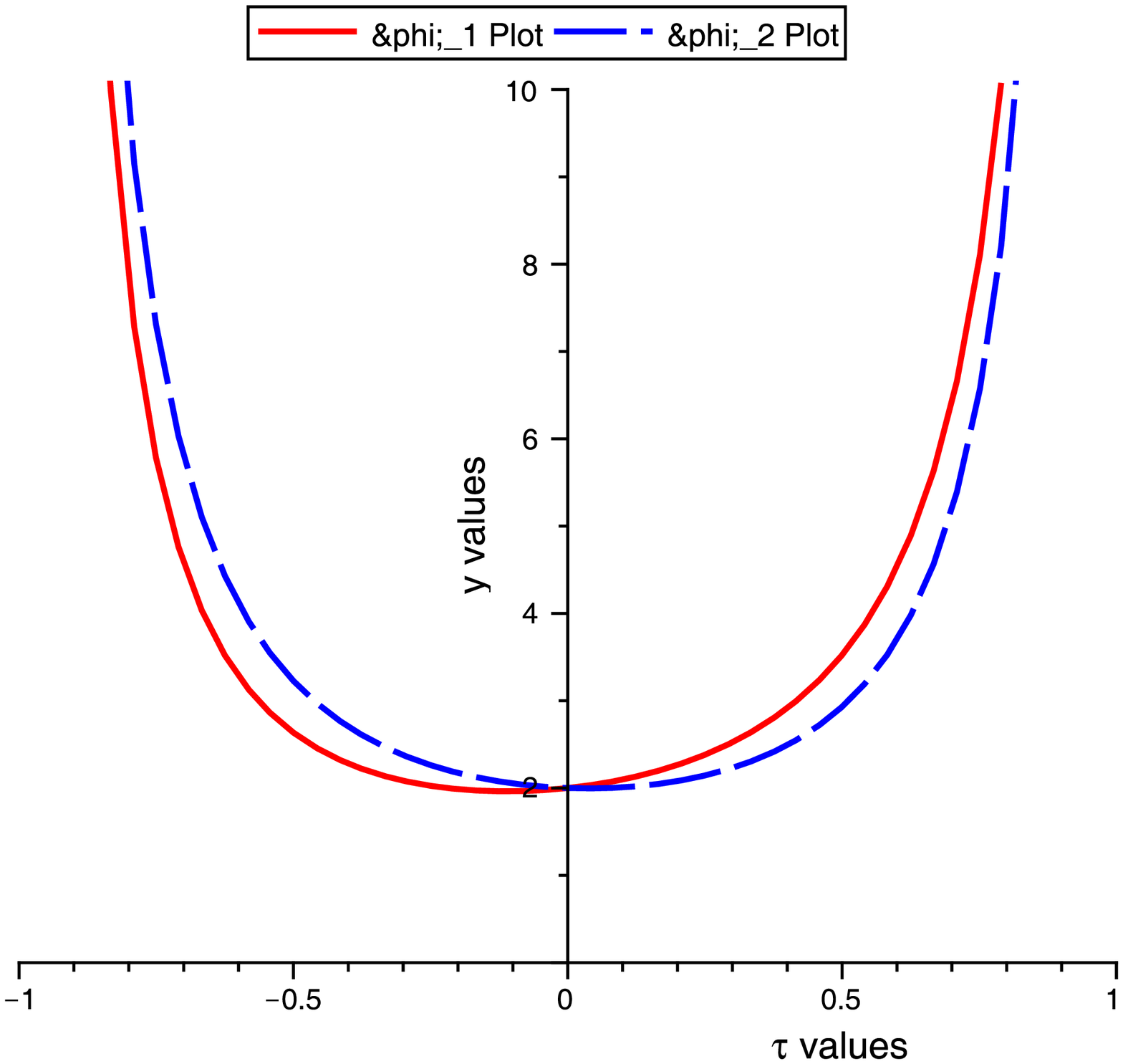}
\\
\caption{ The plots of approximate solutions of Example 1 for M=2.}\label{fig1}
\end{center}
\end{figure}

By the above information the system (\ref{1}) reduces to
\begin{equation}\label{e1-1}
\left\{\begin{array}{ll} \int_{-1}^{1}
\frac{\phi_1(\tau)}{\tau-t}\,d\tau +  \int_{-1}^{1} (\tau-t)
\phi_1(\tau) d\tau + \int_{-1}^{1} t \phi_2(\tau) d\tau = \pi, \ &
-1<t<1
\\ \\
\int_{-1}^{1} \frac{\phi_2(\tau)}{\tau-t}\,d\tau +  \int_{-1}^{1}
\tau \phi_1(\tau) d\tau + \int_{-1}^{1} (\tau+t) \phi_2(\tau) d\tau
= 2 \pi t, \ & -1<t<1
\end{array}\right.
\end{equation}
since the matrices $S=A+B=I_2$ and $D=A-B=-I_2$ are nonsingular,
therefore the solution of system (\ref{e1-1}) exists. The kernels
$K_{1j}(t,\tau)$, $ K_{2j}(t,\tau)$ ($j=1,2$), and the functions
$f_1(t)$, $f_2(t)$ are polynomials of degree at most $1$, so we set
\begin{equation}\label{e1-2}
\phi_j(\tau):=\frac{1}{\sqrt{1-{\tau}^2}}
\left\{\beta_{j0}\,T_0(\tau)+\beta_{j1}\,T_1(\tau)+\beta_{j2}\,T_2(\tau)\right\}\,,
\ \ \ \ j=1,2
\end{equation}
and
\begin{eqnarray*}\begin{array}{lll}
K_{ij}(t,\tau)=\gamma_{ij0}(t) T_0(\tau) + \gamma_{ij1}(t) T_1(\tau)
\,,\ \ \ &i,j=1,2\\
f_i(t)=c_{i0} U_0(t)+c_{i1}U_1(t),\ \ \ &i=1,2& \end{array}
\end{eqnarray*}
where
\begin{eqnarray*}
\begin{array}{cccc}
\gamma_{110}(t)=-\frac{1}{2} U_1(t),   &\gamma_{111}(t)=U_0(t),
&\gamma_{120}(t)=\frac{1}{2}U_1(t),   &\gamma_{121}(t)=0,\\
\gamma_{210}(t)=0, \  &\gamma_{211}(t)=U_0(t), \
&\gamma_{220}(t)=\frac{1}{2} U_1(t),\  &\gamma_{221}(t)=U_0(t),\\
c_{10}(t)=1, \  &c_{11}(t)=0, \
&c_{20}(t)=0,\  &c_{21}(t)=1.\\
\end{array}\end{eqnarray*}
Substituting these expansions into (\ref{e1-1}) and using
(\ref{12})-(\ref{13}), for $\nu=1,$ we obtain
\begin{equation}\label{e1-3}
\left\{\begin{array}{ll}  \beta_{11} U_0(t) +\beta_{12} U_1(t) -
\frac{1}{2} \beta_{10} U_1(t) + \frac{1}{2} \beta_{11} U_0(t) +
\frac{1}{2} \beta_{20} U_1(t) = U_0(t)  \\ \\ \beta_{21} U_0(t)
+\beta_{22} U_1(t) + \frac{1}{2} \beta_{11} U_0(t) + \frac{1}{2}
\beta_{20} U_1(t) + \frac{1}{2} \beta_{21} U_0(t) =
U_1(t)  \\
\end{array}\right.
\end{equation}
Then the linear independency of $\{U_0(t),\ U_1(t)\}$ implies
\begin{eqnarray*}
\left\lbrace
\begin{array}{c} \frac{3}{2} \beta_{11} =1\\
-\frac{1}{2}\beta_{10}+\beta_{12}+\frac{1}{2}\beta_{20} =0\\
\frac{1}{2}\beta_{11}+\frac{3}{2} \beta_{21}
=0\\ \frac{1}{2}\beta_{20}+\beta_{22} =1 .\\
\end{array}\right.
\end{eqnarray*}
A nonunique solution of this system for the arbitrary values of
$\beta_{10}$ and $\beta_{20}$ is given by
\begin{eqnarray*}
\beta_{10},\ \ \ \beta_{11}=\frac{2}{3}, \ \ \
\beta_{12}=\frac{1}{2}\left(\beta_{10}-\beta_{20}\right), \ \ \
\beta_{20}, \ \ \ \beta_{21}=-\frac{2}{9}, \ \ \
\beta_{22}=1-\frac{1}{2}\beta_{20}.
\end{eqnarray*}
For example, if $\beta_{10}=\beta_{20}=2$, then
$\beta_{12}=\beta_{22}=0$ and so we find from (\ref{e1-2})
\begin{eqnarray*}\label{e1-11}
\phi_1(\tau)=\frac{\frac{2}{3} \tau +2}{\sqrt{1-\tau^2}}\,, \ \ \ \
\ \ \phi_2(\tau)=\frac{-\frac{2}{9} \tau +2}{\sqrt{1-\tau^2}}\,,
\end{eqnarray*}
(See figure \ref{fig1} for the behavior of these
solutions).

\textbf{Example $2.$} Solve the problem of Example $1$ in case
$\textbf{3}$.

\begin{figure}
\begin{center}
\includegraphics[width=0.4\textwidth]{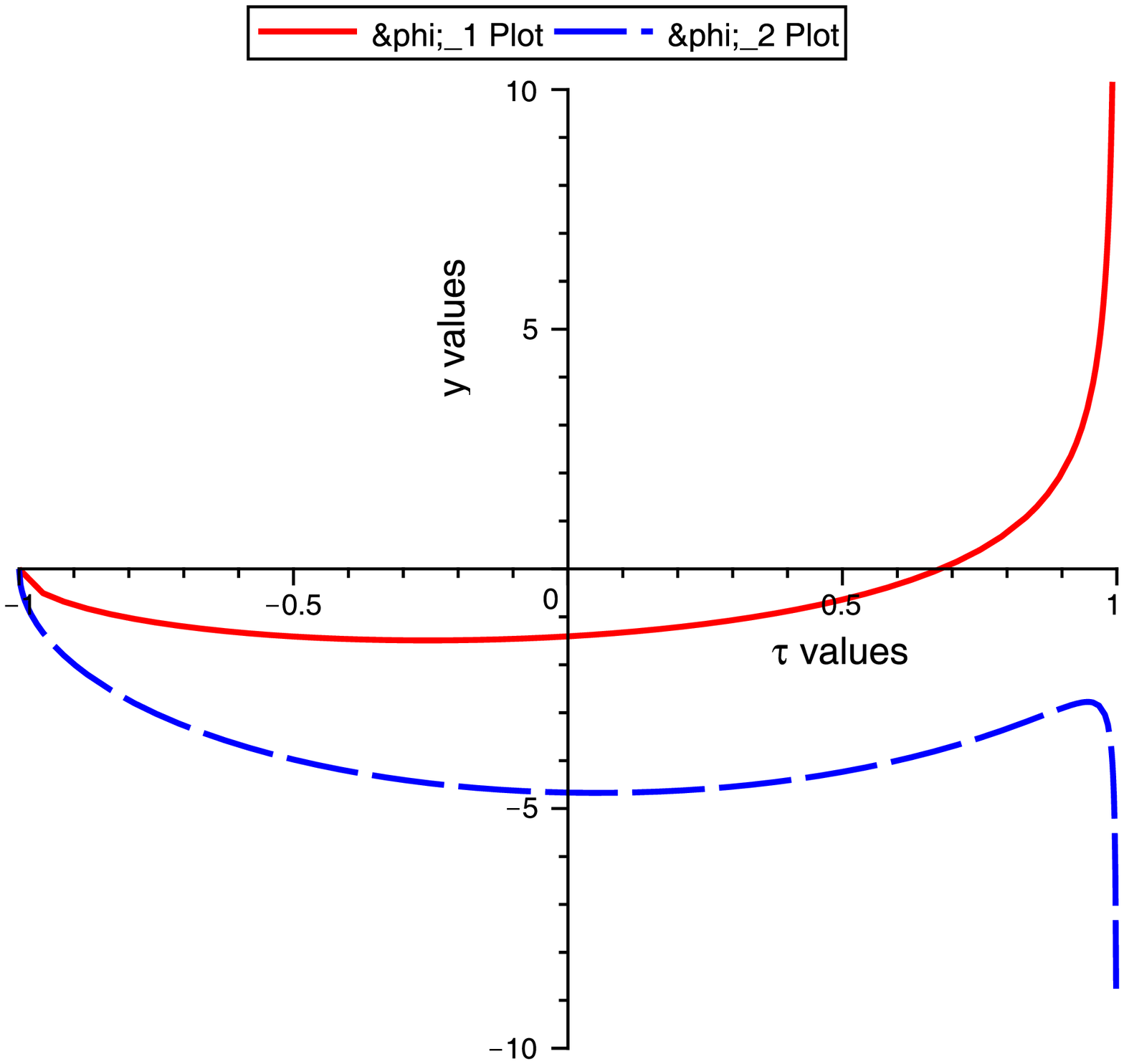}
\\
\caption{The plots of approximate solutions of Example 2 for M=2.}\label{fig2}
\end{center}
\end{figure}

In this case, we set
\begin{equation}\label{e2-2}
\phi_j(\tau):=\sqrt{\frac{1+\tau}{1-\tau}}
\left\{\beta_{j0}\,V_0(\tau)+\beta_{j1}\,V_1(\tau)\right\}\,, \ \ \
\ j=1,2
\end{equation}
and
\begin{eqnarray*}\begin{array}{lll}
K_{ij}(t,\tau)=\gamma_{ij0}(t) V_0(\tau) + \gamma_{ij1}(t) V_1(\tau)
\,,\ \ \ &i,j=1,2,\\
f_i(t)=c_{i0} W_0(t)+c_{i1}W_1(t),\ \ \ &i=1,2,&
\end{array}
\end{eqnarray*}
where
\begin{eqnarray*}
\begin{array}{cccc}
\gamma_{110}(t)=W_0(t)-\frac{1}{2} W_1(t),\ \
&\gamma_{111}(t)=\gamma_{210}(t)=\gamma_{211}(t)=\gamma_{221}(t)=\frac{1}{2}W_0(t),
\\
\gamma_{120}(t)=-\frac{1}{2}W_0(t)+\frac{1}{2}W_1(t),\ &
\gamma_{121}(t)=0,  \ \ \ \ \gamma_{220}(t)=\frac{1}{2}
W_1(t),  \\
c_{10}(t)=1, \ \ \ \ \ \ c_{11}(t)=0, \ & c_{20}(t)=-1, \ \ \ \ \ \
c_{21}(t)=1.
\end{array}\end{eqnarray*}
Substituting these expansions into (\ref{e1-1}) and using
(\ref{12})-(\ref{13}) for $\nu=3$, result
\begin{equation}\label{e2-3}
\left\{\begin{array}{ll}  \beta_{10} W_0(t) +\beta_{11}
W_1(t)+\beta_{10} W_0(t)- \frac{1}{2} \beta_{10} W_1(t),
\\
\ \ \ \ \ \ \ \ \ \ \ \ \ \ \ \ \ \ \ + \frac{1}{2} \beta_{11}
W_0(t) - \frac{1}{2}
\beta_{20} W_0(t)+\frac{1}{2} \beta_{20} W_1(t)= W_0(t)  \\ \\
\beta_{20} W_0(t) +\beta_{21} W_1(t)+\frac{1}{2}\beta_{10} W_0(t)+
\frac{1}{2} \beta_{11} W_0(t) \\
\ \ \ \ \ \ \ \ \ \ \ \ \ \ \ \ \ \ \ \ \ \ \ + \frac{1}{2}
\beta_{20} W_1(t) + \frac{1}{2} \beta_{21} W_0(t)=-W_0(t)+W_1(t)
\end{array}\right.
\end{equation}
and from the linear independency of $\{W_0(t),\ W_1(t)\}$, we get the algebraic system
\begin{eqnarray*}
\left\lbrace
\begin{array}{c}
2\beta_{10}+\frac{1}{2}\beta_{11}-\frac{1}{2}\beta_{20}=1\\
-\frac{1}{2}\beta_{10}+\beta_{11}+\frac{1}{2}\beta_{20}=0\\
\frac{1}{2}\beta_{10}+\frac{1}{2}\beta_{11}+\beta_{20}+\frac{1}{2}\beta_{21}=-1\\
\frac{1}{2}\beta_{20}+\beta_{21} =1,\\
\end{array}\right.
\end{eqnarray*}
which has the unique solution
\begin{eqnarray*}
\beta_{10}=-\frac{10}{27}\,,\,\,\,\beta_{11}=\frac{28}{27}\,,\,\,\,\beta_{20}=-\frac{22}{9}\,,\,\,\,\beta_{21}=\frac{20}{9}\,,
\end{eqnarray*}
and the solutions of (\ref{e1-1}) can be found via (\ref{e2-2}). The graphs of these solutions plotted in Fig. \ref{fig2}.
\\

\begin{figure}
\begin{center}
\includegraphics[width=0.4\textwidth]{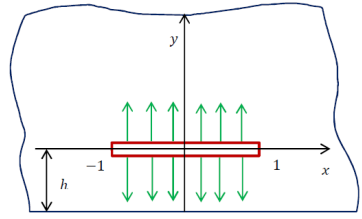}
\\
\caption{Crack parallel to a free boundary}\label{fig3}
\end{center}
\end{figure}

\textbf{Example $3.$} Consider the problem of a half plane
containing a crack parallel to the boundary which illustrated in
Fig. \ref{fig3} and formulated as the system \cite{eg}
\begin{equation}\label{e3-1}
\left\{\begin{array}{ll} \int_{-1}^{1}
\frac{\phi_1(\tau)}{\tau-t}\,d\tau +  \int_{-1}^{1}
\left[K_{11}(t,\tau) \phi_1(\tau) + K_{12}(t,\tau) \phi_2(\tau)
d\tau \right]= 0
\\ \\
\int_{-1}^{1} \frac{\phi_2(\tau)}{\tau-t}\,d\tau +  \int_{-1}^{1}
\left[K_{21}(t,\tau) \phi_1(\tau) + K_{22}(t,\tau) \phi_2(\tau)
d\tau \right]= \pi
\end{array}\right.
\end{equation}
with
\begin{eqnarray*}
K_{11}(t,\tau)=-\frac{\tau-t}{(\tau-t)^2+4h^2}
+\frac{8h^2(\tau-t)}{\left[(\tau-t)^2+4h^2\right]^2}
-\frac{4h^2(\tau-t)\left[12h^2-(\tau-t)^2\right]}{\left[(\tau-t)^2+4h^2\right]^3}\,,\\
K_{12}(t,\tau)=K_{21}(t,\tau)=-\frac{8h^3\left[4h^2-3(\tau-t)^2\right]}{\left[(\tau-t)^2+4h^2\right]^3}\,,\\
K_{22}(t,\tau)=-\frac{\tau-t}{(\tau-t)^2+4h^2}
-\frac{8h^2(\tau-t)}{\left[(\tau-t)^2+4h^2\right]^2}
-\frac{4h^2(\tau-t)\left[12h^2-(\tau-t)^2\right]}{\left[(\tau-t)^2+4h^2\right]^3},\\
\end{eqnarray*}
where $h$ is the distance of crack from the boundary. The  physical
conditions of the problem impose that the relations
\begin{equation}\label{e3-2} \int_{-1}^1
\phi_1(\tau)\,d\tau=0\,, \ \ \ \ \int_{-1}^1 \phi_2(\tau)\,d\tau=0
\end{equation}
and $$\phi_1(t)=\phi_1(-t), \  \ \   \phi_2(t)=-\phi_2(-t)$$
 to be satisfied.
Therefore the unknown functions may be
expressed as
\begin{equation}\label{e3-3}
\phi_1(\tau)\simeq\frac{1}{\sqrt{1-\tau^2}}\,\sum_{j=0}^{M}\,\beta_{1j}\,
T_{2j}(\tau), \ \ \ \
\phi_2(\tau)\simeq\frac{1}{\sqrt{1-\tau^2}}\,\sum_{j=1}^{M}\,\beta_{2j}\,
T_{2j-1}(\tau).
\end{equation}
For $\nu=1$, it follows from the orthogonality condition (\ref{12})
that the second condition in (\ref{e3-2}) satisfies and the first
one gives $\beta_{10}=0$.

By taking $\beta_{11}=\beta_{21}=0$, as arbitrary values,  the
remaining coefficients $\beta_{ij}$ are uniquely determined from the linear algebraic system (\ref{a6}) for each values of $h$ and $M.$ This
leads us to find the functions $\phi_1(\tau)$ and $\phi_2(\tau)$ from
(\ref{e3-3}). \\
The stress intensity factors
\begin{eqnarray*}
k_1=\lim_{\tau\rightarrow 1^{-}} \sqrt{1-\tau^2}\, \phi_2(\tau)
\\
k_2=\lim_{\tau\rightarrow 1^{-}} \sqrt{1-\tau^2}\, \phi_1(\tau)
\end{eqnarray*} and their absolute estimation errors(Est.Err.) reported in table
\ref{tab:1}. For $h=\infty$ and $K_{ij}(t,\tau)=0,$ from Eqs.
(\ref{e3-3}) and (\ref{a5})-(\ref{a6}) the exact solutions of
(\ref{e3-1}) are obtained as
$$ \phi_1(\tau)=0,\ \ \ \ \
\phi_2(\tau)=\frac{\tau}{\sqrt{1-\tau^2}},$$ which give $k_1=1$
and $k_2=0$. This is shown in the last row of Table \ref{tab:1}. The
table shows the rapid convergence of the results even for relatively
small values of $M$.
\begin{footnotesize}
\begin{table}
\caption{Stress intensity factors for the crack parallel to the
boundary}
\label{tab:1}       
\begin{tabular}{llllll}
\hline\noalign{\smallskip}
$h$ & $M$ & $k_1$ & Est.Err. $k_1$ & $k_2$ & Est.Err. $k_2$ \\
\noalign{\smallskip}\hline\noalign{\smallskip}
0.2 &6&4.878800637605022& 6.3e-14& 1.750099102171126&  6.2e-14\\
&7&4.788277537335018& 1.1e-14& 1.727809740547429&   4.1e-15\\
&8&4.760729834685963& 4.8e-14& 1.719782910590219&   1.1e-14\\
\noalign{\smallskip}
0.4&3&2.607272141646415& 4.3e-15& 0.7745787927510580& 1.6e-16\\
&4&2.594500911475041& 7.3e-15& 0.7266641783709941& 5.0e-15\\
&6&2.594423234973139&  4.2e-14& 0.7376171346942053& 3.3e-16\\
\noalign{\smallskip}
0.6&2&1.834057544899021&  1.1e-15& 0.5664257041432605& 4.1e-16\\
&5&1.960455689663461& 6.1e-15& 0.4297949760368867& 1.6e-15\\
\noalign{\smallskip}
0.8&2&1.608371955353828& 1.6e-15& 0.3323260582700188& 2.8e-16\\
&3&1.660617572058080&  8.7e-16& 0.2675691556476836& 4.6e-16\\
\noalign{\smallskip}
1.0&2&1.461157081431933&  2.0e-15& 0.2104682299562445& 1.1e-16\\
&4&1.485914720666516&  2.1e-16& 0.1796691052492212&  1.1e-16\\
\noalign{\smallskip}
1.2&4&1.372176156193755& 5.0e-16& 0.1234414146531335&  0.\\
\noalign{\smallskip}
1.5&4&1.262800608570183& 1.6e-15& 0.07465158121522054& 1.7e-16\\
\noalign{\smallskip}
2.0&3&1.162112249974693& 1.1e-15& 0.03662808437088003& 0.\\
\noalign{\smallskip}
3.0&2&1.077621553329114&  3.3e-16& 0.01274529646673066& 0.\\
\noalign{\smallskip}
10 &2&1.007451045420713&  2.6e-16& 0.00037197952964307& 0.\\
\noalign{\smallskip}
$\infty$ & 1&1&0&0&0\\
\noalign{\smallskip}\hline
\end{tabular}
\end{table}
\end{footnotesize}
\section{Conclusions}
We described a new idea of using Chebyshev polynomials for the numerical solution of system
(\ref{1}). As applications of this idea, we have solved the simple examples $1$ and $2$ as a coupled system of
singular integral equations of the first kind. In example $3$, we studied
a crack problem in solid mechanics and we reported the numerical results (see Table \ref{tab:1}) to show the efficiency
 and rapid convergence of the proposed method for all these kinds of problems.



\end{document}